\documentclass[a4paper,11pt,leqno]{article}
\usepackage[mathscr]{eucal}
\usepackage{amssymb}
\usepackage{latexsym}
\usepackage{theorem}
\theoremstyle{plain}

\newtheorem{Theorem}{Theorem}[section]
\newtheorem{Corollary}{Corollary}[section]
\newtheorem{Lemma}{Lemma}[section]
\newtheorem{Proposition}{Proposition}[section]

\theorembodyfont{\rmfamily}
\newtheorem{Remark}{Remark}[section]

\newtheorem{Definition}{Definition}[section]

\title{Submanifolds with harmonic mean curvature vector 
field in contact $3$-manifolds
\footnote{Colloquium Mathematicum 
\textbf{100} (2004), no.~2,
163--179.
Minor misprints are corrected.}}
\author{Jun-ichi Inoguchi} 
\date{}
\begin{document}
\maketitle
\begin{abstract}
Biharmonic or polyharmonic curves and surfaces in
$3$-dimensional contact manifolds are investigated.
\end{abstract}

\section*{Introduction}

This paper concerns curves and surfaces 
in 3-dimensional contact manifolds,
whose mean curvature vector field is in the
kernel of certain elliptic
differential operators.

First we study submanifolds whose mean curvature
vector field is in the kernel of
Laplacian (submanifolds with
harmonic mean curvature vector fields).

The study of such submanifolds is inspired by 
a conjecture of Bang-yen Chen \cite{Chen}:

\vspace{0.2cm}

{\it Harmonicity of the mean curvature 
vector field implies
harmonicity of the immersion ?}
   
\vspace{0.2cm}

The harmonicity equation $\Delta \mathbb{H}=0$ for
the mean curvature vector field 
$\mathbb{H}$ of an immersed
submanifold
$\mathbf{x}:M^m\to \mathbf{E}^n$ 
in Euclidean $n$-space is equivalent to the
{\it biharmonicity} of the immersion:
$\Delta \Delta \mathbf{x}=0$, since 
$\Delta \mathbf{x}=-m\mathbb{H}$.

A submanifold $\mathbf{x}:M \to\mathbf{E}^n$
is said to be a {\it biharmonic submanifold}
if $\Delta  \mathbb{H}=0$.

In 1985,  Chen proved 
the nonexistence
of proper biharmonic surfaces
in Euclidean 3-space.
The conjecture by Chen is still open.

Some partial and positive answers have
been obtained by several authors
\cite{Defv1}-\cite{Dimitric},
\cite{Garay}-\cite{HV}.

\vspace{0.2cm}

The biharmonicity equation 
is regarded as a special case of the following condition:
$$
\Delta \mathbb{H}=\lambda 
\> \mathbb{H},\ \lambda \in \mathbf{R}.
$$
Namely the mean curvature vector field is 
an eigenfunction of the Laplacian.

The study of Euclidean submanifolds 
with $\Delta \mathbb{H}=\lambda \mathbb{H}$ 
was initiated by
Chen in 1988 (See \cite{Chen}).

It is known that submanifolds
in $\mathbf{E}^n$ satisfying
$\Delta \mathbb{H}=\lambda \mathbb{H}$
are either
biharmonic ($\lambda=0$),
of $1$-type or null $2$-type.
In particular all surfaces in $\mathbf{E}^3$
with $\Delta \mathbb{H}=\lambda \mathbb{H}$
are of constant mean curvature.
Moreover a
surface in $\mathbf{E}^3$
satisfies $\Delta \mathbb{H}=\lambda \mathbb{H}$
if and only if
it is minimal, an open portion of a totally umbilical 
sphere or an open portion of a circular cylinder.

F.~Defever \cite{Defv2}
showed that hypersurfaces satisfying
$\Delta \mathbb{H}=\lambda \mathbb{H}$ are 
of constant mean curvature.
Note that Chen \cite{Chen1994},
\cite{Chen1995} studied spacelike submanifolds
with $\Delta \mathbb{H}
=\lambda \mathbb{H}$ in Minkowski space, 
hyperbolic space or de Sitter space.
M.~Barros and O.~J.~Garay
showed that Hopf cylinders in $S^3$
with $\Delta \mathbb{H}=\lambda \mathbb{H}$
are Hopf cylinders over circles 
in the $2$-sphere $S^2$.
A.~Ferr{\'a}ndez, P.~Lucas and M.~A.~Mero{\~n}o
\cite{FLM}
studied such submanifolds in
anti de Sitter $3$-space $H^3_1$.

\vspace{0.2cm}

In non-constant curvature ambient spaces,
results on biharmonic submanifolds are
very few. 

Recently,
T.~Sasahara \cite{Sasahara2}--\cite{Sasahara3}
studied Legendre surfaces in
the 
Sasakian space form $\mathbf{R}^{5}(-3)$
satisfying $\Delta \mathbb{H}=\lambda \mathbb{H}$.
Moreover
Sasahara introduced the notion of
``$\varphi$-position vector field" 
and ``$\varphi$-mean curvature vector field"
for submanifolds in Sasakian space form
$\mathbf{R}^{2n+1}(-3)$.
Sasahara investigated submanifolds
in $\mathbf{R}^{2n+1}(-3)$
whose $\varphi$-mean curvature vector
field $\mathbb{H}_\varphi$ satisfies
$\Delta \mathbb{H}_\varphi=\lambda \mathbb{H}_\varphi$.
In particular he classified curves and surfaces in
$\mathbf{R}^{3}(-3)$ with 
$\Delta \mathbb{H}_\varphi=\lambda \mathbb{H}_\varphi$.
Since both $\mathbf{R}^{2n+1}(-3)$ and
$S^{2n+1}$ are typical examples of Sasakian space form,
it seems to be interesting to study biharmonic submanifolds
in general Sasakian space forms.

Based on these observations, in the first part of this
paper, we shall study harmonicity of mean curvature
vector fields of  
curves and surfaces in 
3-dimensional Sasakian space forms.
Several results for
3-dimensional sphere $S^3$ due to Spanish 
research group 
(Barros, Garay Ferr{\'a}ndez, Lucas and Mero{\~n}o)
will be generalised
to $3$-dimensional Sasakian space forms.

Next, in the second part,
we shall study another 
``biharmonicity" suggested by 
J.~Eells and J.~H.~Sampson \cite{ES}.
A smooth map $\phi:M\to N$ between Riemannian manifolds
is said to be a biharmonic map (or polyharmonic map
of order $2$)
if its bitension field $\mathscr{T}_{2}(\phi)$ vanishes.
In \cite{CMO}, ``biharmonic" curves and surfaces in $S^3$
 are classified.
We shall classify Legendre curves and Hopf cylinders 
in $3$-dimensional
Sasakian space forms, which are biharmonic in this sense.

In particular we shall show the existence of
non-minimal biharmonic
Hopf cylinders in Sasakian space forms of holomorphic
sectional curvature greater than $1$ (Berger spheres).

The author would like to thank Dr.~Cezar
 Dumitru Oniciuc (University ``AL. I. Cuza ") 
 and 
 Dr.~Tooru Sasahara (Hokkaido University)
 for their useful comments.

\section*{Part I}
\section{Preliminaries}
\subsection{Contact manifolds}\label{contact}
We begin by recalling fundamental 
ingredients of contact Riemannian geometry
from \cite{Blair}.

Let $M$ be a $(2n+1)$-manifold. 
A one form $\eta$ is called a {\it
contact form} on $M$ if $(d\eta)^{n} \wedge \eta \not=0$. 
A $(2n+1)$-manifold
$M$ together with a contact form is called a {\it contact manifold}. 
The {\it contact distribution} $D$ of $(M,\eta)$ is defined by
$$
D=\left \{ X \in TM \ \vert \
\eta(X)=0
\right \}.
$$

On a contact manifold $(M,\eta)$, there exists a unique vector
field $\xi$ such that
$$
\eta (\xi)=1,\ \  d \eta (\xi,\cdot)=0.
$$
This
vector field $\xi$ is called the {\it Reeb vector field} or {\it
characteristic vector field} of $(M,\eta)$.

Moreover there exists an endomorphism field $\varphi$ and a
Riemannian metric $g$ on $M$ such that
\begin{equation}\label{almostcontact}
\varphi^2=-I+\eta \otimes \xi, \ \eta(\xi)=1,
\end{equation}
\begin{equation}\label{metric}
g(\varphi X,\varphi Y)=g(X,Y)-\eta(X)\eta(Y),
\ \ g(\xi,\cdot)=\eta,
\end{equation}
\begin{equation}\label{associatedmetric}
d\eta(X,Y)=2g(X,\varphi Y)
\end{equation}
for all vector fields $X,\ Y$ on $M$. 
On an almost contact manifold $(M,\eta;\xi,\varphi)$, there exists
a Riemannian metric $g$ satisfying (\ref{metric}). 
Such a metric $g$ is
called an {\it compatible metric} of $M$.
A contact manifold $(M,\eta)$ together with 
structure tensors $(\xi,\varphi,g)$ is called 
a {\it contact Riemannian manifold}.

\begin{Proposition}
Let $(M,\eta,\xi,\varphi,g)$ be a 
contact Riemannian manifold.
\newline
\noindent
Then $M$ $\xi$ is a Killing vector field
if and only if
\begin{equation}\label{K-contact}
\nabla_{X}\xi=-\varphi X,\ \
X \in \mathfrak{X}(M).
\end{equation}
Here $\nabla$ is the Levi-Civita connection of $(M,g)$.
\end{Proposition}

\begin{Definition}
A contact Riemannian manifold 
$(M,\eta,\xi,\varphi,g)$
is 
said to be a {\it Sasaki manifold} if 
\begin{equation}\label{Sasakian}
(\nabla_{X}\varphi)Y
=g(X,Y)\xi-\eta(Y)X, \ \ X,Y \in
\mathfrak{X}(M).
\end{equation}
\end{Definition}

Note that on a Sasaki manifold, 
$\xi$ is a Killing vector field.

\vspace{0.2cm}

Let $(M,\eta;\xi,\varphi,g)$ be a
contact Riemannian manifold. 
A tangent plane at a point of $M$ 
is said to be a {\it holomorphic
plane} if it is invariant under $\varphi$. The sectional curvature
of a holomorphic plane is called {\it holomorphic sectional
curvature}. If the sectional curvature function of $M$ is constant
on all holomorphic planes in $TM$, then $M$ is said to be of {\it
constant holomorphic sectional curvature}.
Complete and connected Sasaki manifolds of constant holomorphic
sectional curvature are called {\it Sasakian space forms}.
Let us denote by $R$ the Riemannian curvature tensor
field of the metric $g$
which is defined by
$$
R(X,Y):=\nabla_{X}\nabla_{Y}-
\nabla_{Y}\nabla_{X}-
\nabla_{[X,Y]},\ \ 
X,Y \in \mathfrak{X}(M).
$$
When $(M,\eta;\xi,\varphi,g)$ is a
Sasakian space form of constant 
holomorphic sectional curvature $c$, then $R$
is 
described by the
following formula:

\begin{eqnarray*}
R(X,Y)Z & = & \frac{c+3}{4}
\{
g(Y,Z)X-g(Z,X)Y \}\\
& & +\frac{c-1}{4}\>\{
\eta(Z)\eta(X)Y
-\eta(Y)\eta(Z)X \\
& &+g(Z,X)\eta(Y)\xi
-g(Y,Z)\eta(X)\xi \\
& &-g(Y,\varphi Z)
\varphi X-g(Z,\varphi X)\varphi Y+
2g(X,\varphi Y)\varphi Z\
\>\}.
\end{eqnarray*}

Note that even if the 
holomorphic sectional curvature is negative, 
a Sasakian space form is {\it not} negatively curved. 
In fact, the sectional curvature of plane sections
containing $\xi$ is $1$ on any Sasaki manifold.

\vspace{0.2cm}


It is known that every $3$-dimensional 
Sasakian space form is realised as a 
Lie group together with a left invariant Sasaki structure. 
More precisely the following is known ({\it cf}. \cite{BTV}):
\begin{Proposition}
Simply connected $3$-dimensional 
Sasakian space form of 
constant holomorphic sectional curvature is isomorphic to 
\begin{enumerate}
\item special unitary group $\mathrm{SU}(2);$
\item Heisenberg group $\mathbf{R}^{3}(-3);$
\item the universal covering group of the special linear group
$\mathrm{SL}_{2}\mathbf{R}$
\end{enumerate}
\noindent together with canonical left invariant Sasaki structure.
In particular simply connected
Sasakian space form of constant holomorphic
sectional curvature $1$ is the 
$\mathrm{SU}(2)$ with biinvariant metric of constant
curvature $1$ (hence isometric to the
unit $3$-sphere $S^3$). 
\end{Proposition}

\subsection{Boothby-Wang fibration}

Let $(M^{2n+1},\eta;\xi,\varphi,g)$ be a
contact Riemannian
manifold. Then
$M$ is said to be {\it regular} if $\xi$ generates a one-parameter
group $K$ of isometries on $M$, such that the action of $K$
on $M$ is simply transitive.
Note that if $M$ is regular, then both
$\varphi$ and $\eta$ are automatically $K$-invariant, {\it i.e},
$\pounds_{\xi}\varphi=0$ and $\pounds_{\xi}\eta=0$. The Killing
vector field $\xi$ induces a regular one-dimensional Riemannian
foliation on $M$. We denote by ${\overline M}:=M/K$ the orbit space
(the space of all leaves) of a regular contact Riemannian manifold $M$
under the $K$-action.

Let $\bar{X}_{\bar p}$ be a tangent 
vector of the orbit space
$\overline{M}$ at ${\bar p}=\pi(p)$. 
Then there exists a tangent vector 
${\bar X}^{*}_p$ of $M$ at $p$ which is
orthogonal to $\xi$  such that
$\pi_{*p}{\bar X}^{*}_p=\bar{X}_{\bar p}$.
The tangent vector ${\bar X}^{*}_p$ 
is called the
{\it horizontal lift} 
of $\bar{X}_{\bar p}$ to $M$
at $p$.
The horizontal lift operation
$*:{\bar X}_{\bar p} \mapsto 
{\bar X}^{*}_{p}$
is naturally extended to vector fields. 

The contact structure on $M$ 
induces an almost Hermitian structure 
on the orbit space ${\overline M}$:
\begin{equation}\label{basecpxstr}
J{\bar X}=\pi_{*}(\varphi {\bar X}^{*}),\ 
\bar{X} \in \mathfrak{X}(\bar{M}).
\end{equation}

Let us denote by $\bar{\nabla}$ the Levi-Civita connection
of $\bar{M}$. Then, 
by using the fundamental equations
for Riemannian submersions due to 
O'Neill \cite{ON},
we have the following results.

\begin{Proposition}{\rm(\cite{Og})}
Let $M$ be a regular contact Riemannian manifold.
Then for any $\bar{X},\bar{Y}
\in \mathfrak{X}(\bar{M}):$
\begin{equation}\label{submersion}
\nabla_{{\bar X}^*}\bar{Y}^{*}=
(\bar{\nabla}_{\bar X}{\bar Y})^{*}
-g(\bar{X}^{*},\varphi \bar{Y}^{*})\xi.
\end{equation}
\end{Proposition}

\begin{Proposition}{\rm(\cite{Og})}
Sasakian space forms are regular Sasaki manifolds.
The orbit space of a Sasakian space form of 
constant holomorphic sectional curvature $c$ is a 
complex space form of constant holomorphic 
sectional curvature $c+3$.
\end{Proposition}
W.~M.~Boothby and H.~C.~Wang \cite{BW} 
proved that if $M$
is a compact regular contact manifold, 
then the natural projection
$\pi:M \rightarrow {\bar M}$ defines a principal circle bundle
over a symplectic manifold ${\bar M}$ and the symplectic form
$\Omega$ of ${\overline M}$ 
determines an integral cocycle. Furthermore
the contact form $\eta$ gives a connection form of this circle
bundle and satisfies $\pi^{*}\Omega=d\eta$. The fibering $\pi:M
\rightarrow {\bar M}$ is called the {\it Boothby-Wang fibering}
of a regular compact contact manifold $M$.
Based on this result, we call the fibering $\pi:M \to \bar{M}$
of a regular contact Riemannian manifold $M$, the 
``Boothby-Wang fibering"
of $M$ even if $M$ is noncompact.

The unit sphere  $S^{2n+1}$ is a typical example of regular compact
Sasaki manifold. For $S^{2n+1}$, the Boothby-Wang fibering
coincides with the {\it Hopf fibering} $S^{2n+1}\rightarrow
\mathbb{C}P^n$.

In $3$-dimensional case, the Boothby-Wang fibering 
of Sasakian space forms  
have the following matrix group models \cite{BTV}:
\begin{eqnarray*}
\pi:\mathrm{SU}(2)\to S^2(c)&=&\mathrm{SU}(2)/\mathrm{U}(1)
,\\ 
\pi: \mathbf{R}^{3}(-3)\to \mathbf{C}&=&
\mathbf{R}^{3}(-3)/\mathbf{R},\\
 \pi: \mathrm{SL}_{2}\mathbf{R} \to H^{2}( c)
&=&\mathrm{SL}_{2}\mathbf{R}/\mathrm{SO}(2).
\end{eqnarray*}
Here $S^2(c)$ and $H^2 (c)$ 
are sphere and 
hyperbolic space of curvature $c$, 
respectively.

\subsection{Hopf cylinders}\label{Hopfcylinders}

Now we shall restrict our attention to $3$-dimensional regular
contact Riemannian manifold $M$.

Let ${\bar \gamma}$ be a curve parameterized by arc length
in ${\overline M}$ with curvature ${\bar \kappa}$. Taking the
inverse image $S_{\bar \gamma}:=\pi^{-1}\{{\bar \gamma}\}$ 
of ${\bar \gamma}$
in $M^3$.

Here we compute the fundamental quantities
of $S_{\bar \gamma}$.

Let us denote by $\bar{P}=({\bar \mathbf p}_1,{\bar \mathbf p}_2)$
the Frenet frame field of $\bar{\gamma}$. 
By using the complex structure $J$ of
${\overline M}^2$,
${\bar \mathbf{p}}_2$
is given by
$$
{\bar \mathbf{p}}_2=J
{\bar \mathbf{p}}_1
$$
Then the Frenet-Serret formula
of $\bar{\gamma}$ is given by
$$
\bar{\nabla}_{\bar{\gamma}^\prime}P=P
\left(
\begin{array}{cc}
0 &-\bar{\kappa}  \\
\bar{\kappa} &0
\end{array}
\right).
$$

Here the function $\bar{\kappa}$
is the (signed) curvature of $\bar{\gamma}$.

Let $\mathbf{t}=(\bar{\mathbf{p}}_1)^{*}$
the horizontal lift of 
$\bar{\mathbf{p}}_1$ with respect to the
Boothby-Wang fibering.
Then $(\mathbf{t},\xi)$ gives an orthonormal
frame field of $S$.
We choose a unit normal vector
field $\mathbf{n}$ by
$\mathbf{n}
=(\bar{\mathbf p}_2)^{*}$. 
Since ${\bar {\mathbf p}}_2$ is defined by
${\bar \mathbf{p}}_2=J{\bar \mathbf{p}}_1$,
$\mathbf{n}=\varphi \> \mathbf{t}$.
In fact, 
$$
(\bar{\mathbf p}_2)^{*}=
(J\bar{\mathbf p}_1)^{*}=\varphi
(\bar{\mathbf p}_1)^{*}=\varphi
\> \mathbf{t}.
$$

Let us denote by $\nabla^S$ the Levi-Civita connection
of $S$. 
The {\it second fundamental form}
$I\!I$
derived from $\mathbf{n}$
is defined by the
{\it Gau{\ss} formula}:
\begin{equation}\label{}
\nabla_{X}Y=\nabla^{S}_{X}Y+
I\!I(X,Y)\mathbf{n},\ \
X, Y \in \mathfrak{X}(S).
\end{equation}

By using (\ref{submersion}),
$$
\nabla_{\mathbf t}\>\mathbf{t}=
(\bar{\nabla}_{\bar{\mathbf p}_{1}}\bar{\mathbf p}_{1})^{*}
-g(\mathbf{t},\varphi\>\mathbf{t})\xi
=(\bar{\kappa}\circ \pi)\mathbf{n}.
$$
Hence $\nabla^{S}_{\mathbf t}{\mathbf t}=0$.
Since $\xi$ is Killing, we have
$\nabla^{S}_{\mathbf t}\xi=\nabla^{S}_{\xi}\xi=0$.
Thus $S_{\bar \gamma}$ is flat.
The second fundamental form 
$I\!I$ is described as
$$
I\!I(\mathbf{t},\mathbf{t})=\bar{\kappa} \circ \pi,
\ \
I\!I(\mathbf{t},\xi)=-1,\ \
I\!I(\xi,\xi)=0.
$$
The mean curvature is $H=(\bar{\kappa} \circ \pi)/2$
and the mean curvature vector field $\mathbb{H}$ is
$\mathbb{H}=H\> \mathbf{n}$.

In case $M=S^3$, $S_{\bar \gamma}$ is called the {\it Hopf cylinder}.
In particular if ${\bar \gamma}$ is closed, 
then $S_{\bar \gamma}$ is a flat torus in $S^3$
and 
called the {\it Hopf torus} over ${\bar \gamma}$
(H.~B.~Lawson, {\it cf}. \cite{Kitagawa}, \cite{Pin}).
The Hopf torus over a geodesic in 
$S^2(4)$ coincides with the Clifford minimal torus.
We call the flat surface $S_{\bar \gamma}$ in a regular 
contact Riemannian manifold
$M$ a {\it Hopf cylinder} over the curve $\bar{\gamma}$ in 
$\overline{M}$.

\subsection{Curves in Riemannian $3$-manifolds}\label{CurvesSasaki}
Let $(M,g)$ be a Riemannian manifold
and $\gamma=\gamma(s):I\to M$
a curve parametrised by the arclength parameter
in $M$. We regard $\gamma$ as a 
1-dimensional Riemannian 
manifold with respect to the metric induced by
$g$.

We recall the following definition ({\it cf}. \cite{BB}). 
\begin{Definition}
If $\gamma(s)$ is a unit speed curve 
in a Riemannian $3$-manifold $(M^3,g)$,
we say that $\gamma$ is a {\it Frenet curve} if there exists
an orthonomal frame field $P=(\mathbf{p}_1,\mathbf{p}_2,
\mathbf{p}_3)$ along $\gamma$ and two nonnegative functions
$\kappa$ and $\tau$ such that
$P$ satisfies the following {\it Frenet-Serret formula}:
$$
\nabla_{\gamma^\prime}P=P
\left(
\begin{array}{ccc}
0 & -\kappa & 0 \\
\kappa & 0 & -\tau \\
0 & \tau & 0
\end{array}
\right), \ \
\mathbf{p}_1=\gamma^{\prime}(s).
$$ 
The functions $\kappa$ and 
$\tau$ are called
the {\it curvature} and {\it torsion}
of $\gamma$ respectively.
\end{Definition}

Geodesics can be regarded as Frenet curves with 
$\kappa=0$.  
A curve with constant curvature and zero torsion 
is called a ({\it Riemannian}) {\it circle}.
A helix is a curve whose curvature and torsion are constants.
Riemannian circles are regarded as degenerate helices.
Helices, which are not circles, are frequently
called {\it proper helices}.

Note that, in general ambient space $(M^3,g)$,
geodesics may have non-vanishing torsion. In fact,
as we shall see later, Legendre geodesics in a
Sasakian 3-manifold have constant torsion $1$.

The Frenet-Serret formula of $\gamma$
implies that the mean curvature vector field $\mathbb{H}$ 
of a Frenet curve $\gamma$ is
given by
$$
\mathbb{H}=\nabla_{\gamma^\prime}\gamma^\prime
=\kappa \mathbf{p}_2.
$$

Let us denote by $\Delta$ the Laplace operator
acting on the space $\Gamma(\gamma^{*}TM)$
of all smooth sections of the vector bundle:
$$
\gamma^{*}TM:=\bigcup_{s\in I}T_{\gamma(s)}M
$$
over $I$.
Then $\Delta$ is given explicitly by
$$
\Delta=-\nabla_{\gamma^\prime}
\nabla_{\gamma^\prime}.
$$

\begin{Lemma}
The mean curvature vector field
 $\mathbb{H}$ of a Frenet curve $\gamma$
is harmonic  in 
$\gamma^{*}TM$ 
{\rm(}$\Delta \mathbb{H}=0${\rm)} if 
and only if
$$
\nabla_{\gamma^\prime}
\nabla_{\gamma^\prime}
\nabla_{\gamma^\prime}
\gamma^\prime=0.
$$
\end{Lemma}

\vspace{0.2cm}

When $M$ is the 
Euclidean space $\mathbf{E}^m$,
a curve $\gamma$ satisfies $\Delta \mathbb{H}=0$
if and only if $\gamma$ is
{\it biharmonic}, {\it i.e}., $\Delta \Delta \gamma=0$ since
$\Delta \gamma=-\mathbb{H}$.

The following general result is essentially obtained in
\cite{FLM}.

\begin{Theorem}\label{harmonicH}
Let $\gamma$ be a Frenet curve in a Riemannian
$3$-manifold $(M,g)$.
Then $\gamma$ satisfies
$\Delta \mathbb{H}=\lambda \mathbb{H}$
in $\gamma^{*}TM$ if and only if  
$\gamma$ is a 
geodesic $(\lambda=0)$
or a helix satisfying
$\lambda=\kappa^{2}+\tau^{2}$.  
\end{Theorem}
{\it Proof.}
Let $I$ be an open interval and
$\gamma=\gamma(s):I \to M$ be a curve 
parametrised by
the arclength parameter $s$ with Frenet frame field 
$P=(\mathbf{p}_1,\mathbf{p}_2,\mathbf{p}_3)$.
Direct computation shows that

\begin{equation}\label{nabla-H}
\nabla_{\gamma^\prime}\mathbb{H}=
-\kappa^{2}\mathbf{p}_1+
\kappa^\prime \>
\mathbf{p}_{2}+\kappa \tau \mathbf{p}_3.
\end{equation}

Let us compute the Laplacian 
of $\mathbb{H}$:
$$
-\Delta \mathbb{H}=
\nabla_{\gamma^\prime}\nabla_{\gamma^\prime}\mathbb{H}=
-3\kappa\kappa^\prime \> \mathbf{p}_1
+(\kappa^{\prime \prime}-\kappa^3-\kappa \tau^2)
\mathbf{p}_2+
(2\kappa^{\prime}\tau+\kappa \tau^{\prime})\mathbf{p}_3.
$$
Hence $\Delta \mathbb{H}=\lambda \mathbb{H}$ if and only if
$$
\kappa \> \tau^{\prime}=0,\ \ \kappa^{3}+\kappa \> \tau^{2}=\lambda \kappa.
$$
These formulae imply that 
$\gamma$ is a geodesic or a helix 
satisfying $\lambda=\kappa^{2}+\tau^{2}$.

Conversely every geodesic 
satisfies $\Delta \mathbb{H}=0$. 
Helices satisfy $\Delta \mathbb{H}=\lambda
\mathbb{H}$ with
$\lambda=\kappa^{2}+\tau^{2}$.
$\Box$

\begin{Corollary}{\rm(\cite{Dimitric})}
Let $\gamma$ be a curve in Euclidean $3$-space $\mathbf{E}^3$.
Then $\gamma$ is biharmonic if and only if $\gamma$ is a 
straight line.
\end{Corollary}

On the contrary, in indefinite 
semi-Euclidean space, there exist nongeodesic biharmonic curves.
Chen and Ishikawa \cite{ChenIshikawa}
classified 
biharmonic spacelike curves
in $\mathbf{E}^m_{\nu}$.
(See also \cite{Inoguchi-L}).

\subsection{Curves with normal-harmonic mean curvature}

The results in the preceding subsection say that to characterise 
curves which are non geodesics
we need to use another
differential operator for our purpose.

In this subsection we use the {\it normal} Laplacian.

Let $\gamma:I \to M$ be a Frenet curve in 
an oriented Riemannian $3$-manifold
$M$ parametrised by the arclength.
Denote by 
$P=(\mathbf{p}_1,\mathbf{p}_2,\mathbf{p}_3)$ 
the Frenet frame field of $\gamma$ as before.
Then the {\it normal bundle} $T^{\perp}\gamma$ of the
curve $\gamma$ is given by
$$
T^{\perp}\gamma=\bigcup_{s\in I}T^{\perp}_{s}\gamma,\
T^{\perp}_{s}\gamma=\mathbf{R}\> \mathbf{p}_{2}(s)\oplus 
\mathbf{R}\> \mathbf{p}_{3}(s).
$$
The {\it normal connection}
$\nabla^{\perp}$
is a connection of $T^{\perp}\gamma$
defined by
$$
\nabla_{\gamma^\prime}^{\perp}X=\mathrm{normal}\
\mathrm{component}\
\mathrm{of}\
\nabla_{\gamma^\prime}X
$$
for any section $X$ of the normal
bundle $T^{\perp}\gamma$.

By using the Frenet frame field, 
$\nabla^{\perp}$ can be represented as
$$
\nabla^{\perp}_{\gamma^\prime}X=
\nabla_{\gamma^\prime}X-
g(\nabla_{\gamma^\prime}X,\mathbf{p}_1)
\mathbf{p}_1.
$$

Let us denote by $\Delta^{\perp}$
the Laplace operator 
acting on
the space $\Gamma(T^{\perp}\gamma)$
of all smooth sections of the normal
bundle $T^{\perp}\gamma$.
The operator $\Delta^{\perp}$ is called the
{\it normal Laplacian}
of $\gamma$ in $M$.
The normal Laplacian $\Delta^\perp$ is
given by
$$
\Delta^{\perp}X=-\nabla^{\perp}_{\gamma^\prime}
\nabla^{\perp}_{\gamma^\prime}X,\ \ X \in \Gamma(T^{\perp}\gamma).
$$

\vspace{0.2cm}

Now we compute $\Delta^{\perp}\mathbb{H}$.
From (\ref{nabla-H}), we have
$$
\nabla^{\perp}_{\gamma^\prime}\mathbb{H}=
\kappa^\prime \>
\mathbf{p}_{2}+\kappa \tau \mathbf{p}_3.
$$
From this equation, we get
$$
-\Delta^{\perp}\mathbb{H}=
(\kappa^{\prime \prime}-\kappa \tau^{2})\mathbf{p}_{2}+
(2\kappa^{\prime}\tau+\kappa \tau^\prime)\mathbf{p}_{3}.
$$
\begin{Theorem}{\rm({\it cf.\/} \cite{FLM})}\label{normaleigen}
A curve $\gamma$ satisfies 
$\Delta^{\perp}\mathbb{H}=
\lambda \mathbb{H}$
if and only if
$$ 
\kappa^{\prime \prime}-\kappa \tau^{2}=
-\lambda \kappa
,\ \
2\kappa^{\prime}\tau+\kappa \tau^{\prime}=0.
$$
\end{Theorem}

\begin{Corollary}\label{normalharmonic}
A curve $\gamma$ satisfies $\Delta^{\perp}\mathbb{H}=0$
if and only if
$$ 
\kappa^{\prime \prime}-\kappa \tau^{2}=0,\ \
2\kappa^{\prime}\tau+\kappa \tau^{\prime}=0.
$$
\end{Corollary}

We shall apply these general results for curves in 
Sasakian $3$-manifolds in the next section.
Note that Barros and Garay classified curves, which satisfy
$\Delta^{\perp}\mathbb{H}=\lambda \mathbb{H}$ in space forms
\cite{BG2},\cite{BG3}.

\section{Curves and surfaces in $3$-dimensional \\
Sasaki manifolds}

\subsection{Curves in $3$-dimensional Sasaki manifolds}

Now let $M^3=(M,\eta,\xi,\varphi,g)$ be a 
contact Riemannian $3$-manifold
with an associated metric $g$.
A curve $\gamma=\gamma(s):I \to M$
parametrised by the arclength parameter is 
said to be a {\it Legendre curve} if $\gamma$ is tangent to
the contact distribution $D$ of $M$.
It is obvious that $\gamma$ is Legendre if and only if
$\eta(\gamma^{\prime})=0$.

Let $\gamma$ be a Legendre curve in $M^3$.
Then we can take a Frenet frame field 
$P=(\mathbf{p}_1,\mathbf{p}_2,\mathbf{p}_3)$
so that $\mathbf{p}_1=\gamma^\prime$ and
$\mathbf{p}_3=\xi$. (See \cite{BB}).

Now we assume that $M$ is a Sasaki manifold.
Then by definition,
the Frenet-Serret formula of $\gamma$ is given explicitly by
$$
\nabla_{\gamma^\prime}P=P
\left(
\begin{array}{ccc}
0 & -\kappa & 0 \\
\kappa & 0 & -1 \\
0 & 1 & 0
\end{array}
\right).
$$
Namely every Legendre curve has constant torsion $1$
\cite{BB}.

\vspace{0.2cm}

Now we investigate curves with harmonic or normal-harmonic
mean curvature vector field in Sasakian 3-manifolds.

The following two results are direct
consequence of Theorem \ref{harmonicH}
and Theorem \ref{normaleigen},
respectively.

\begin{Corollary}
Let $\gamma$ be a Legendre curve in $3$-dimensional
Sasaki manifold.
Then $\gamma$ satisfies
$\Delta \mathbb{H}=\lambda \mathbb{H}$
in $\gamma^{*}TM$ 
if and only if $\gamma$ is a 
Legendre geodesic $(\lambda=0)$ 
or a Ledendre helix satisfying 
$\lambda=\kappa^{2}+1$ $(\lambda\not=0)$.
\end{Corollary}

\begin{Remark}
Sasaki manifolds together with
compatible Lorentz metric are called
{\it Sasakian spacetimes}
(\cite{Duggal},\cite{Takahashi}). On Sasakian spacetimes,
the Reeb vector fields are timelike.
Every $3$-dimensional Sasakian spacetime contains
proper biharmonic Legendre curves.
In fact, in a $3$-dimensional Sasakian spacetime
biharmonic Legendre curves are 
Legendre geodesics or
Legendre helices with curvature $1$.
 ({\it cf.\/} \cite{Inoguchi-L}). 
\end{Remark}

\begin{Proposition}
Let $\gamma$ be a Legendre curve in a Sasakian $3$-manifold.
Then $\Delta^{\perp}\mathbb{H}=\lambda \mathbb{H}$ 
if and only if
$\gamma$ is a Legendre geodesic $(\lambda=0)$
or a Legendre helix with constant nonzero curvature
$(\lambda\not=0)$.
In the latter case, $\lambda=1$.
\end{Proposition}

\subsection{Biharmonic Hopf cylinders}\label{biharmHopf}

In this section we study harmonicity and
normal-harmonicity of the mean curvature
of Hopf cylinders.

Let $M^3$ be a regular Sasaki manifold
with Boothby-Wang fibration 
$\pi:M \to \bar{M}$.

Take a curve $\bar{\gamma}=\bar{\gamma}(s)$
parametrised by the arclength $s$ in the base
space form $\bar{M}$.
Let us denote by $S=S_{\bar \gamma}$ the Hopf cylinder
of $\bar{\gamma}$. 
(See Section \ref{Hopfcylinders})

Let $\mathbf{t}=(\bar{\mathbf{p}}_1)^{*}$ be
the horizontal lift of 
$\bar{\mathbf{p}}_1$ with respect to the
Boothby-Wang fibering.
Then $(\mathbf{t},\xi)$ gives an orthonormal
frame field of $M$.
The unit normal vector field
$\mathbf{n}$ is the horizontal lift
of $\bar{\mathbf p}_2$. 
Note that
$\mathbf{n}=\varphi \mathbf{t}$.

The mean curvature vector field $\mathbb{H}$ of $S$ is 
$\mathbb{H}=H\> \mathbf{n}=(\bar{\kappa}\circ \pi)\mathbf{n}/2$.

\vspace{0.2cm}

Now we study harmonicity and normal-harmonicity
of $\mathbb{H}$. 
Denote by $\iota$ the inclusion map
of $S$ into $M$.
Then the Laplace operator
$\Delta$ acting on the space
$\Gamma(\iota^{*}TM)$ and
the
normal Laplacian $\Delta^\perp$
of $S$ are given by
$$
\Delta=
-\left(
\nabla_{\mathbf t}\nabla_{\mathbf t}+
\nabla_{\xi}\nabla_{\xi}
\right),
\ \
\Delta^\perp=
-\left(
\nabla^{\perp}_{\mathbf t}\nabla^{\perp}_{\mathbf t}+
\nabla^{\perp}_{\xi}\nabla^{\perp}_{\xi}
\right),
$$
respectively. 
Direct computation shows that
$$
\nabla_{\mathbf t}\mathbb{H}=-2H^2\mathbf{t}+H^{\prime}\mathbf{n}+H\xi,
 \ \
\nabla^{\perp}_{\mathbf t}\mathbb{H}=H^{\prime}\mathbf{n},
\ \
\nabla_{\xi}\> \mathbb{H}=H \>\mathbf{t},
\ \
\nabla^{\perp}_{\xi}\> \mathbb{H}=0,
$$
$$
\nabla_{\xi}\nabla_{\xi}\mathbb{H}=-H \mathbf{n}.
$$
Thus we get
$$
-\Delta \mathbb{H}=-6HH^{\prime}\mathbf{t}+
(H^{\prime \prime}-4H^{3}-2H)\mathbf{n}+
2H^{\prime}\xi, 
$$
$$
-\Delta^{\perp} \mathbb{H}=
H^{\prime \prime}\mathbf{n}.
$$

\begin{Theorem}
A Hopf cylinder $S_{\bar \gamma}$ in a
$3$-dimensional regular Sasaki manifold satisfies
$\Delta \mathbb{H}=\lambda \mathbb{H}$
in $\iota^{*}TM$ if and only if
$\bar{\gamma}$ is a geodesic $(\lambda=0)$ or a Riemannian circle
$(\lambda \not=0)$.
In case that $\lambda\not=0$, the eigenvalue $\lambda$ is
$\lambda=4H^2+2> 2$.
\end{Theorem}

\begin{Remark}
Every Hopf cylinder in a $3$-dimensional
regular Sasaki manifold is {\it anti invariant}.
Sasahara showed that
an anti invariant surface in $\mathbf{R}^{3}(-3)$ satisfies
$\Delta \mathbb{H}=\lambda \mathbb{H},\ \lambda\not=0$
if and only if it is a Hopf cylinder over a circle with 
$\lambda>2$.
See Proposition 11 in \cite{Sasahara2}.
\end{Remark}

\begin{Lemma}
A Hopf cylinder $S_{\bar \gamma}$ satisfies
$\Delta^{\perp} \mathbb{H}=\lambda \mathbb{H}$ if and only if
$\gamma$ is defined by one of the following
natural equations{\rm :}
\begin{enumerate}
\item
$\bar{\kappa}(s)=as+b, \ a, b\in \mathbf{R}, \ \lambda=0;$
\item
$\bar{\kappa}(s)=a \cos (\sqrt{\lambda}s)+
b \sin (\sqrt{\lambda}s),\ \lambda >0;$
\item
$\bar{\kappa}(s)=a \exp (\sqrt{-\lambda}s)+
b \exp (-\sqrt{-\lambda}s),\ \lambda <0.$
\end{enumerate}
\end{Lemma}
{\it Proof.}
The Hopf cylinder $S_{\bar \gamma}$ satisfies
$\Delta^{\perp} \mathbb{H}=\lambda \mathbb{H}$
if and only if $\bar{\gamma}$ satisfies
$$
{\bar \kappa}^{\prime \prime}+
\lambda {\bar \kappa}=0.
$$
Thus the result follows. $\Box$

\begin{Theorem}\label{Hopfnormalharmonic}
A Hopf cylinder $S_{\bar \gamma}$ satisfies
$\Delta^{\perp} \mathbb{H}=0$ if and only if
$\bar{\gamma}$ is one of the following{\rm :}
\begin{enumerate}
\item
a geodesic{\rm ;}
\item
a Riemannian circle or{\rm ;}
\item
a Riemannian clothoid
{\rm(} Cornu spiral {\rm)}.
\end{enumerate}
Here a Riemannian clothoid is a curve in $\bar{M}^2$
whose curvature is a linear function of the arclength.
\end{Theorem}

\begin{Remark}
On curves in Riemannian 2-space forms,
the following result is 
obtained \cite{FLM}:
\end{Remark}

\begin{Theorem}
Let $\bar{\gamma}$ be a curve in 
Riemannian $2$-manifold
$\bar{{M}}^2$.
To avoid the confusion,
let us denote by $\Delta^{\perp}_{\bar \gamma}$
and $\mathbb{H}_{\bar \gamma}$
the normal Laplacian of $\bar{\gamma}$ 
and the mean curvature vector
in 
$\bar{{M}}^2$ respectively.
Then $\Delta_{\bar \gamma}^{\perp}\mathbb{H}_{\bar \gamma}
=\lambda \mathbb{H}_{\bar \gamma}$ if and only if
\begin{enumerate}
\item $\bar{\gamma}$ is a geodesic,
Riemannian circle or a Riemannian clothoid{\rm ;}
\item
$\bar{\kappa}(s)=a \cos (\sqrt{\lambda}s)+
b \sin (\sqrt{\lambda}s),\ \lambda >0${\rm ;}
\item
$\bar{\kappa}(s)=a \exp (\sqrt{-\lambda}s)+
b \exp (-\sqrt{-\lambda}s),\ \lambda <0$.
\end{enumerate}
\end{Theorem}
\begin{Corollary}
Let $M$ be a $3$-dimensional regular Sasaki manifold 
\newline
\noindent
with
Boothby-Wang fibering $\pi:M \to \bar{M}$.
Let $\bar{\gamma}$ be a curve in 
$\bar{M}$. Then the Hopf cylinder
$S=S_{\bar \gamma}$ satisfies
$\Delta^{\perp}\mathbb{H}=\lambda \mathbb{H}$
if and only if $\bar{\gamma}$ satisfies
$\Delta_{\bar \gamma}^{\perp}\mathbb{H}_{\bar \gamma}
=\lambda \mathbb{H}_{\bar \gamma}$.
\end{Corollary}
Theorem \ref{Hopfnormalharmonic} is a generalisation of
a result obtained by Barros and Garay \cite{BG}.
In fact, if we choose $M^3=S^3$ then we obtain the 
following.

\begin{Theorem}\label{HopfnormalharmonicS3}
{\rm (\cite{BG})}
A Hopf cylinder $S_{\bar \gamma}$ in the
unit $3$-sphere $S^3$ satisfies
$\Delta^{\perp} \mathbb{H}=0$ if and only if
$\gamma$ is one of the following{\rm :}
\begin{enumerate}
\item
a geodesic, 
\item
a Riemannian circle or
\item
a Riemannian clothoid.
\end{enumerate}
Here a Riemannian clothoid is a curve in the
$2$-sphere $S^2(1/2)$ of radius $1/2$
whose curvature is a linear function of the arclength.
\end{Theorem}

Riemannian clothoids 
are called 
``Cornu spirals" in
\cite{BG}.

\section*{Part II}
\section{Polyharmonic maps}

Let $(M^m,g)$ and $(N^n,h)$ be Riemannian manifolds
and $\phi:M \to N$ a smooth map.
The {\it tension field}
$\mathscr{T}(\phi)$ is a section
of the vector bundle $\phi^{*}(TN)$
defined by
$$
\mathscr{T}(\phi):=\mathrm{tr}(\nabla d \phi).
$$
A smooth map $\phi$ is said to be 
a {\it harmonic map} if its tension field vanishes.
It is well known that $\phi$ is harmonic if and
only if $\phi$ is a critical 
point of the {\it energy}:
$$
E(\phi)=\int
\frac{1}{2}|d\phi |^{2}dv_{g}
$$
over every compact supported region of $M$.

Now let $\phi$ be a harmonic map.
Then the Hessian $\mathcal{H}_{\phi}$
of the energy is given by the
following second variation formula:
$$
\mathcal{H}_{\phi}(V,W)=\int
h(\mathcal{J}_{\phi}(V),W)dv_{g},\ \ 
V,W \in \Gamma(\phi^{*}TN).
$$
Here the operator $\mathcal{J}_{\phi}$
is the {\it Jacobi operator}
of the harmonic map $\phi$ defined by
$$
\mathcal{J}_{\phi}(V):=
\bar{\Delta}_{\phi}V-\mathcal{R}_{\phi}(V),\
\ V \in \Gamma(\phi^{*}TN),
$$
$$
\bar{\Delta}_{\phi}:=-\{\sum_{i=1}^{m}(
\nabla^{\phi}_{e_i}\nabla^{\phi}_{e_i}-\nabla^{\phi}_
{\nabla_{e_i}e_i}\},\ \
\mathcal{R}_{\phi}(V)=\sum_{i=1}^{m}R^{N}(V,d\phi (e_i))
d\phi(e_i).
$$
Here $\nabla^\phi,\ R^N$ and $\{e_i\}$ denote the
induced connection of $\phi^{*}TN$,
 curvature tensor of $N$
and a local orthonormal
frame field of $M$,
respectively.

For general theory of harmonic maps and their
Jacobi operators, we refer to \cite{EL} and \cite{Ur}.

J.~Eells and J.~H.~Sampson suggested to 
study {\it polyharmonic maps}
(See \cite{ES} and \cite{EL}, p.~77 (8.7)).
Let $\phi:M \to N$ be a smooth map as before.
Then $\phi$ is said to be
a polyharmonic map of order $k$
if it is an extremal
of the functional:
$$
E_{k}(\phi)=\int |(d+d^{*})^{k}\phi |^{2}dv_{g}.
$$
Here $d^{*}$ is the codifferential operator.
In particular, if $k=2$, we have
$$
E_{2}(\phi)=\int |\mathscr{T}(\phi)|^{2} dv_{g}.
$$

The Euler-Lagrange equation
of the functional $E_{2}$
was computed by Caddeo and Oproiu
(See \cite{CMO}, p.~867) and G.~Y.~Jiang \cite{Jiang1}--\cite{Jiang2},
independently.
The Euler-Lagrange equation of $E_{2}$ is
$$
\mathscr{T}_{2}(\phi):=-\mathcal{J}_{\phi}(\mathscr{T}(\phi))=0.
$$

\begin{Remark}
Let $\phi:M\to N$ be an isometric immersion.
Then its tension field is $m\mathbb{H}$.
Thus the functional $E_{2}$ is given by
$$
E_{2}(\phi)=m^{2}\int |\mathbb{H}|^{2}dv_{g}.
$$
In case that $M$ is $2$-dimensional, $E_{2}(\phi)$
the total mean curvature of $M$  
up to constant multiple. See \cite{Chenbook}, Section 5.3. 
\end{Remark}

In particular, if $N=\mathbf{E}^n$ and $\phi$
an isometric immersion, then
$$
\mathscr{T}_{2}(\phi)=-\Delta_{M} \Delta_{M} \phi,
$$ 
since $\Delta_{M} \phi=m \mathbb{H}$.
Here $\Delta_{M}$ is the Laplacian of $(M,g)$.
Thus the polyharmonicity (of order $2$) for an isometric immersion
into Euclidean space is equivalent to the
biharmonicity in the sense of Chen.
On this reason, polyharmonic maps of order $2$ are frequently called
{\it biharmonic maps} (or $2$-harmonic maps)
\cite{CMO}, \cite{Jiang1},
\cite{Jiang2}, \cite{Oniciuc}.

Obviously, the notion of $p$-harmonic map
in the sense of \cite{ELrep2}, p.~397 is different from
that of polyharmonic map of order $p$. 

Hereafter we call polyharmonic maps of order $2$ 
by the name ``polyharmonic maps" in short.

Caddeo, Montaldo and Oniciuc 
classified polyharmonic curves 
in \newline
\noindent
$3$-dimensional
Riemannian space forms. More precisely they showed the
following two results.

\begin{Theorem}{\rm (\cite{CMO})}
Let $N$ be a $3$-dimensional
Riemannian space form of nonpositive curvature.
Then all the polyharmonic curves 
are geodesics.
\end{Theorem}

Next for the study 
of polyharmonic curves
in positively curved space forms,
we may assume that $N^3$ is the unit $3$-sphere.

\begin{Theorem}{\rm (\cite{CMO})}
Let $\gamma:I \to S^3$ be a polyharmonic curve
parametrised by the arclength.
Then $0\leq \kappa \leq 1$ and
$\gamma$ is one of the following{\rm:}
\begin{enumerate}
\item $\kappa$ is a geodesic {\rm(}$\kappa=0${\rm ) }.
\item If $k=1$ then $\gamma$ is a
Riemannian circle of curvature $1;$
\item If $0 < \kappa <1$ then $\gamma$ 
is a geodesic of the
Clifford minimal torus of $S^3$.
\end{enumerate}
\end{Theorem}

The preceding theorem implies the following result:

\begin{Corollary}\label{Legendrebiharmonicsphere}
Let $\gamma:I \to S^3$ be a Legendre curve parametrised by
the arclength. Then $\gamma$ is polyharmonic if and only if
$\gamma$ is a Legendre geodesic.
\end{Corollary}

In fact, curves in the the latter two classes 
can not be Legendre. (Recall that every Legendre curve
has constant torsion $1$). 

Now we study polyharmonic Legendre curves 
in contact Riemannian $3$-manifolds.

Let $M^3$ be a contact Riemannian $3$-manifold
and $\gamma:I \to M$ a Frenet curve framed by
$(\mathbf{p}_{1},\mathbf{p}_{2},\mathbf{p}_{3})$.
Then direct computation shows that
$$
\mathscr{T}_{2}(\gamma)=
-3\kappa\kappa^{\prime}\mathbf{p}_{1}
+(\kappa^{\prime \prime}-\kappa^{3}-\kappa \tau^2)\mathbf{p}_{2}
+(2\kappa^{\prime}\tau+\kappa \tau^{\prime})\mathbf{p}_{3}
+\kappa R(\mathbf{p}_{2},\mathbf{p}_{1})\mathbf{p}_{1}.
$$

Now assume that $M$ is a Sasakian space form
of constant holomorphic sectional curvature $c$ then
\begin{eqnarray*}
R(\mathbf{p}_{2},\mathbf{p}_{1})\mathbf{p}_{1}&=&
\frac{c+3}{4}\mathbf{p}_{2} \\
&+&\frac{c-1}{4}
\{\eta(\mathbf{p}_{2})
\eta(\mathbf{p}_{1})\mathbf{p}_{1} \\
&-&\eta(\mathbf{p}_{1})^{2}\mathbf{p}_{2}-
\eta(\mathbf{p}_{2})\xi
+3g(\mathbf{p}_{2},\varphi 
\mathbf{p}_{1})\varphi \mathbf{p}_{1}
\}.
\end{eqnarray*}

In particular, if $\gamma$ is Legendre,
then $R({\mathbf p}_{2},{\mathbf p}_{1}){\mathbf p}_{1}
=c\>\mathbf{p}_{2}$.
Thus a Legendre curve
$\gamma$ in $M$ is polyharmonic if and only if
$$
\kappa=\mathrm{constant},\
\kappa^{3}-(c-1)\kappa=0,\ \tau=1.
$$
If we look for 
nongeodesic polyharmonic Legendre curves, we obtain
$$
\kappa=\mathrm{constant}, \
\kappa^{2}=c-1,\ 
\tau=1.
$$
Thus we obtain the following result
which is a generalisation of Corollary 
\ref{Legendrebiharmonicsphere}.

\begin{Theorem}
Let $M^3(c)$ be a Sasakian space form
of constant holomorphic sectional
curvature $c$ and $\gamma:I \to M$ a
polyharmonic Legendre curve
parametrised by the arclength.

\begin{enumerate}
\item If $c\leq 1$, 
then $\gamma$ is a Legendre geodesic{\rm ;}
\item If $c>1$, 
then $\gamma$ is a Legendre geodesic or a
Legendre helix 
of curvature $\sqrt{c-1}$.
\end{enumerate}
\end{Theorem}

Let $\phi:M \to N$ be an isometric immersion.
Then $\phi$ is a critical point of the
volume functional if and only if 
$\phi$ is minimal.
The {\it Jacobi operator} $\mathscr{J}$
of a minimal immersion $\phi$ (with respect to
the volume functional)
is appeared in the second variation formula
of the volume and given by \cite{Simons}
$$
\mathscr{J}V=\Delta^{\perp}V-\mathscr{S}V+\mathscr{R}(V),\
\ V \in \Gamma(T^{\perp}M).
$$
Here the operators $\mathscr{S}$ and
$\mathscr{R}$ are defined by
$$
h(\mathscr{S}V,W)=\mathrm{tr}(\mathcal{A}_V \circ \mathcal{A}_{W}),\
\
\mathscr{R}(V)=\sum_{i=1}^{m}(R^{N}(d\phi(e_{i}),V)
d \phi (e_{i}))^{\perp}.
$$
Here $\mathcal{A}_V$ denotes the Weingarten operator
with respect to $V$.

Arroyo, Barros and Garay studied
submanifolds in $S^3$ whose
mean curvature vector fields are eigen-section
of the Jacobi operator with respect to the
{\it volume functional}
\cite{ABG}, \cite{BG2},
\cite{BG3}. Such study for surfaces in 
$5$-dimensional Sasakian space forms can be found in
\cite{Sasahara2}.

It seems to be interesting to study similar problems for
submanifolds in space forms or Sasakian space forms
with respect to the energy functional.

In \cite{CMO}, all the polyharmonic 
surfaces in $S^3$ are classified.
More precisely, the only non-minimal polyharmonic surfaces are
totally umbilical $2$-spheres.

Based on this result, we would like to propose the
following problem:

\vspace{0.2cm}

{\it Are there non-minimal 
and non totally umbilical polyharmonic 
submanifolds in homogeneous Riemannian manifolds ?}

\vspace{0.2cm}

To close this paper, 
we study polyharmonic Hopf cylinders in 
\newline
\noindent
$3$-dimensional Sasakian space forms.
Moreover we show the existence of
non-minimal and 
non totally umbilical
polyharmonic surfaces  
in Sasakian space forms.

First we recall the following result which is a
consequence of the main result in \cite{CMO}:

\begin{Proposition}
There are no non minimal polyharmonic Hopf cylinders in
the unit $3$-sphere $S^3$.
\end{Proposition}

Now we generalise this result to Sasakian space forms.

Let $S=S_{\bar \gamma}$ be a Hopf cylinder and 
$\iota:S \subset M^{3}(c)$
its inclusion map into a Sasakian space form $M^{3}(c)$.
Then the bitension field $\mathscr{T}_{2}(\iota)$ is given by
$$
\mathscr{T}_{2}(\iota)
=-\mathcal{J}_{\iota}(\mathscr{T}(\iota))
=-2\mathcal{J}_{\iota}(\mathbb{H}).
$$
We use the orthonormal frame field $\{{\mathbf t},\xi\}$ as before.
Then since $S$ is flat, we have
$$
\bar{\Delta}_{\iota}\mathbb{H}=\Delta \mathbb{H},
\ \
\mathcal{R}(\mathbb{H})=H(
R(\mathbf{n},\mathbf{t})\mathbf{t}+
R(\mathbf{n},\xi)\xi
).
$$
Using the curvature formula of Sasakian space form,
we get
$$
\mathcal{R}(\mathbb{H})=(c+1)H\> \mathbf{n}.
$$
Hence 
$$
\mathcal{J}_{\iota}(\mathbb{H})=6HH^{\prime}
\mathbf{t}-
(H^{\prime \prime}-4H^{3}+(c-1)H)\mathbf{n}
-2H^{\prime}\xi.
$$

Thus $\mathcal{J}_{\iota}(\mathbb{H})=\lambda \> \mathbb{H}$ if and only if
$$
H^{\prime}=0,\ 4H^{3}=(c-1+\lambda)H
$$
and hence $H=0$ or $\lambda=4H^{2}+1-c, \ H\not=0$.
\begin{Theorem}
Let $S$ be a Hopf cylinder 
in a Sasakian space form $M^{3}(c )$.
Then $S$ satisfies 
$\mathcal{J}_{\iota}(\mathbb{H})
=\lambda \> \mathbb{H}$ 
if and only if the base curve of $S$ is a Riemannian circle or a geodesic.
In case that the base curve is not a geodesic, then $\lambda=4H^{2}+1-c$.
\end{Theorem}
\begin{Corollary}
Let $\iota:S_{\bar \gamma}\to M^{3}(c)$ be a 
polyharmonic Hopf cylinder in a Sasakian space form.
\begin{enumerate}
\item
If $c\leq 1$ then $\bar{\gamma}$ is a geodesic{\rm;}
\item
If $c>1$ then $\bar{\gamma}$ 
is a geodesic or a Riemannian circle of curvature
$\bar{\kappa}=\sqrt{c-1}$.
\end{enumerate}
In particular, 
there exist  
nonminimal polyharmonic Hopf cylinders in
\newline
\noindent
Sasakian space forms of holomorphic sectional
curvature greater than $1$.
\end{Corollary}


Department of Mathematics Education,
Faculty of Education,

Utsunomiya University, Minemachi 350,
Utsunomiya, 321-8505,
Japan

{\it E-mail address}: {\tt inoguchi@cc.utsunomiya-u.ac.jp}

\vspace{0.3cm}

Current Address:

\vspace{0.2cm}

Department of Mathematical Sciences,
Faculty of Science,

Yamagata University, Kojirakawa 1-4-12,
Yamagata, 990-8560,
Japan

{\it E-mail address}: {\tt inoguchi@sci.kj.yamagata-u.ac.jp}

\end{document}